\documentclass[12pt]{article}
\usepackage{amssymb, amsmath}
\begin{document}
\renewcommand{\theequation}{\arabic{section}.\arabic{equation}}
\newtheorem{theorem}{Theorem}[section]
\newtheorem{lemma}{Lemma}[section]
\newtheorem{pro}{Proposition}[section]
\newtheorem{cor}{Corollary}[section]
\newcommand{\n}{\nonumber}
\newcommand{\tv}{\tilde{v}}
\newcommand{\tw}{\tilde{\omega}}
\renewcommand{\t}{\theta}
\newcommand{\w}{\omega}
\renewcommand{\b}{\dot{B}^0_{\infty ,1}}
\newcommand{\bn}{\|_{\dot{B}^0_{\infty ,1}}}
\newcommand{\e}{\varepsilon}
\renewcommand{\a}{\alpha}
\renewcommand{\l}{\lambda}
\newcommand{\vare}{\varepsilon}
\newcommand{\s}{\sigma}
\renewcommand{\o}{\omega}
\renewcommand{\O}{\Omega}
\newcommand{\bb}{\begin{equation}}
\newcommand{\ee}{\end{equation}}
\newcommand{\bq}{\begin{eqnarray}}
\newcommand{\eq}{\end{eqnarray}}
\newcommand{\bqn}{\begin{eqnarray*}}
\newcommand{\eqn}{\end{eqnarray*}}
\title{Nonexistence of asymptotically self-similar singularities in the  Euler and the
 Navier-Stokes equations}
\author{Dongho Chae\thanks{This work was supported partially by the KOSEF Grant no.
R01-2005-000-10077-0. \newline {\bf Keywords}: Euler equations,
Navier-Stokes equations, self-similar singularity\newline
}\\
Department of Mathematics\\
              Sungkyunkwan University\\
              Suwon 440-746, Korea\\
  e-mail: {\it chae@skku.edu}}
 \date{}
  \maketitle
\begin{abstract}
In this paper we rule out the possibility of asymptotically
self-similar singularities for  both of the 3D Euler and  the 3D
Navier-Stokes equations. The notion means that the local in time
classical solutions of the equations develop self-similar profiles
as $t$ goes to the possible time of singularity $T$. For the Euler
equations we consider the case where the vorticity converges to
the corresponding self-similar voriticity profile in the sense of
the critical Besov space norm, $\dot{B}^0_{1, \infty}(\Bbb R^3)$.
For the Navier-Stokes equations the convergence of the velocity to
the self-similar singularity is in $L^q(B(z,r))$  for some $q\in
[2, \infty)$, where the ball of radius $r$ is shrinking toward a
possible singularity point $z$ at the order of $\sqrt{T-t}$ as $t$
approaches to $T$. In the $L^q (\Bbb R^3)$ convergence case with
$q\in [3, \infty)$ we present a simple alternative proof of the
similar result in \cite{hou}.
\end{abstract}

\section{Introduction}
 \setcounter{equation}{0}
The problems of global in time regularity/finite time singularity
in the 3D Euler equations and the 3D Navier-Stokes equations are
among the most important and at the same time the most challenging
open problems in the mathematical fluid mechanics(see e.g.
\cite{con1, con2, maj1, maj2}). For rather general introduction to
 the mathematical theories of the Euler and the
Navier-Stokes equations we refer \cite{maj1, che,tem,con4,lad}.
Although there are many partial progresses for the Euler
equations(e.g. \cite{bea, con3}) and for the Navier-Stokes
equations(e.g. \cite{caf, ser, oky, pro, esc, lin}), the solutions
to the problem still look too far to be seen. On the other hand,
in many of the nonlinear partial differential equations where  the
finite time singularity is searched for, one of the most popular
scenario  to check is by the self-similar ansatz, consistent with
the scaling properties of the equations. For the 3D Navier-Stokes
this type of possibility leading to a self-similar singularity was
first considered by Leray in \cite{ler}, and its nonexistence was
proved in \cite{nec}, and the result was later refined by the
authors in \cite{tsa, mil}. For the 3D Euler equations similar
nonexistence result has  been recently obtained by the author of
this article in \cite{cha1}.
 More refined notion of `asymptotically self-similar singularity' is
 considered by the authors in \cite{gig2}, in the context of the
 nonlinear scalar heat equation, and also by physicists including
 the authors of \cite{gre,pel} in the context of 3D Euler equations.
 It means that the local
in time smooth solution evolves into a self-similar profile as the
possible singularity time is approached. The meaning of it will be
more clear in the statements of Theorem 1.2, Theorem 1.4 and
Theorem 1.5 below. For the 3D Navier-Stokes equations the similar
notion was considered rather indirectly by Hou and Li(\cite{hou}),
and they obtained the nonexistence result, assuming that the
convergence of the local in time smooth solution to the
self-similar profile occurs in the $L^q(\Bbb R^3)$ sense in terms
of the self-similar variables with $3<q<\infty$. In this paper,
concerning the case of Navier-Stokes equations, we prove
nonexistence result for the asymptotically self-similar
singularity, for which our convergence of the local classical
solution to the self-similar singularity is in $L^q(B(z,r))$ sense
 with more flexible range of $q\in [2, \infty)$, where the ball of radius
$r\propto\sqrt{T-t}$ is shrinking to a possible singularity point
$z$ as $t$ approaches to $T$. This could be regarded as a
localized version of the corresponding result of \cite{hou}.
Moreover, in the global convergence case in $L^q(\Bbb R^3)$, we
present here an alternative  simple proof of \cite{hou},
 using a classical result on the blow-up rate estimate due to Leray in
 \cite{ler}.
 We apply our argument also to prove nonexistence of asymptotically
  self-similar singularity for
the 3D Euler equations. Here, the use of critical Besov space
$\dot{B}^0_{\infty, 1}(\Bbb R^3)$, which is slightly more regular
than $L^\infty (\Bbb R^3)$, is crucial to obtain our results. The
remaining part of introduction will be divided into two
subsections, and we provide some preliminaries on the Besov
spaces, and state the main theorems for the Euler and the
Navier-Stokes equations respectively.

 \subsection{The Euler Equations}
The system of Euler equations for the homogeneous incompressible
fluid flows in $\Bbb R^3$ is the following.
 \[
\mathrm{ (E)}
 \left\{ \aligned
 &\frac{\partial v}{\partial t} +(v\cdot \nabla )v =-\nabla p,
 \quad (x,t)\in {\Bbb R^3}\times (0, \infty) \\
 &\textrm{div }\, v =0 , \quad (x,t)\in {\Bbb R^3}\times (0,
 \infty)\\
  &v(x,0)=v_0 (x), \quad x\in \Bbb R^3
  \endaligned
  \right.
  \]
where $v=(v_1, v_2, v_3 )$, $v_j =v_j (x, t)$, $j=1,2,3$, is the
velocity of the flow, $p=p(x,t)$ is the scalar pressure, and $v_0 $
is the given initial velocity, satisfying div $v_0 =0$.  We
introduce briefly the homogeneous Besov spaces,
$\dot{B}^s_{\infty,1} (\Bbb R^n)$ and its inhomogeneous
 counterpart,
 $B^s_{\infty,1}(\Bbb R^n)$. For more details on the Besov spaces
 we refer \cite{tri}.
 Given $f\in \mathcal{S}$, the Schwartz class of rapidly deceasing
 functions in $\Bbb R^n$, its Fourier transform $\hat{f}$ is defined by
 $$
 \mathcal{F}(f)= \hat{f} (\xi)=\frac{1}{(2\pi )^{n/2}}\int_{\Bbb R^n} e^{-ix\cdot \xi }
 f(x)dx.
  $$
 We consider  $\varphi \in \mathcal{S}$ satisfying the following three conditions:
\begin{itemize}
 \item[(i)] $\textrm{Supp}\, \hat{\varphi} \subset
 \{\xi \in {\Bbb R}^n  \,| \,\frac12 \leq |\xi|\leq
 2\}$,
 \item[(ii)] $\hat{\varphi} (\xi)\geq C >0 $ if $\frac23 <|\xi|<\frac32$,
 \item[(iii)] $
 \sum_{j\in \mathbb{Z}}  \hat{\varphi}_j (\xi )=1$,
  where $\hat{\varphi_j } =\hat{\varphi } (2^{-j} \xi
 )$.
 \end{itemize}
Construction  of such sequence of functions $\{ \varphi_j \}_{j\in
\Bbb Z}$ is well-known(see e.g. \cite{tri}). Given $s\in \Bbb R$,
the norm of the homogeneous Besov space $\dot{B}^s_{\infty,1} (\Bbb
R^n)$ is now defined by
 $$ f\in \dot{B}^s_{\infty ,
 1}(\Bbb R^n) \Longleftrightarrow \|f\|_{\dot{B}^s_{\infty,1}}:=
 \sum_{j\in \Bbb Z}2^{sj}\|\varphi_j *
 f\|_{L^\infty} < \infty,
 $$
 where $*$ denotes the convolution,
 $(f*g)(x)=\int_{\Bbb R^n} f(x-y)g(y)dy$. The norm $\|\cdot
 \|_{\dot{B}^s_{\infty,1}}$  is actually defined up to addition of
 polynomials(namely,
  if $f_1-f_2 $ is a polynomial, then both of $f_1$ and $f_2$ give the same norm),
 and the space  $\dot{B}^s_{\infty,1} (\Bbb R^n)$ is defined as the
 quotient space of  a class of functions with finite norm, $\|\cdot\|_{\dot{B}^s_{\infty,1}}$,
 divided by the space of polynomials in $\Bbb R^n$.
 Let us set $\hat{\Phi}(\xi )= \sum_{k\leq -1} \varphi _k (\xi)$ if $\xi\neq 0$, and
  define $\hat{\Phi} (0)=1 $.  Then the inhomogeneous
 Besov space is defined by
$$ f\in B^s_{\infty ,
 1}(\Bbb R^n) \Longleftrightarrow \|f\|_{B^s_{\infty ,
 1}}:= \|\Phi* f\|_{L^\infty}+\sum_{j \geq 0} 2^{sj}\|\varphi_j *
 f\|_{L^\infty} < \infty.
 $$
Note that the condition (iii)  and above definitions imply
immediately that both of the spaces
 $ \dot{B}^0_{\infty ,1}(\Bbb R^n)$ and $B^0_{\infty , 1}(\Bbb R^n)$ are continuously embedded into $L^\infty (\Bbb
R^n)$.  The space $B^0_{\infty , 1}(\Bbb R^n)$, in particular,
actually can be embedded into the class of continuous bounded
functions, thus having slightly better regularity than $L^\infty
(\Bbb R^n)$, but containing as a subspace the H\"{o}lder space
$C^{0,\gamma} (\Bbb R^n) $, for any $\gamma>0$. We begin with
statement of a new type of continuation principle for local in time
classical solutions of the Euler system.

\begin{theorem}
Let $v\in C([0, T);B^1_{\infty, 1}(\Bbb R^3))$ be a classical
solution to the 3D Euler equations. There exists an absolute
constant $\eta >0$ such that if
 \bb\label{th11}
 \inf_{0\leq t<T} (T-t) \|\o (t)\|_{\dot{B}^0_{\infty, 1} } <\eta
 ,
 \ee
 the, $v\in C([0, T+\delta); B^1_{\infty, 1}(\Bbb R^3))$ for some $\delta
 >0$.
\end{theorem}
 \noindent{ \textsf{Remark 1.1}} The proof of the local
existence for $v_0 \in B^1_{\infty, 1} (\Bbb R^3)$  is
  implied in the proofs of the main theorems in \cite{cha3, cha4}(see also \cite{vis}), and
 explicitly written in \cite{pak}.  The above theorem implies that
 if $T_*$ is the first time of singularity, then  we have
 the lower estimate of the blow-up rate,
 \bb\label{blow}
 \|\o (t)\|_{\dot{B}^0_{\infty, 1} } \geq \frac{C}{T_*-t}\quad
 \forall t\in [0, T_*)
 \ee
 for an absolute constant $C$. The estimate (\ref{blow}) was actually
 derived previously by a different argument in \cite{cha4}.
 We observe that (\ref{blow}) is consistent both with  the
 Beale-Kato-Majda criterion(\cite{bea}) and
 with Kerr's numerical calculation in \cite{ker} respectively.\\

\noindent{ \textsf{Remark 1.2}} The above continuation principle
 for a local solutions in $B^1_{\infty, 1}$ has  obvious
application to the  solutions belonging to more conventional
function spaces, due to the embeddings,
$$ H^m (\Bbb R^3 )\hookrightarrow C^{1, \gamma} (\Bbb R^3)\hookrightarrow
 B^1_{\infty, 1}(\Bbb R^3)$$
 for $m>5/2$ and  $\gamma = m-3/2$.
For example the local solution $v\in C([0, T); H^m (\Bbb R^3 ))$
can be continued to be $v\in C([0, T+\delta);  H^m (\Bbb R^3 ))$
for some $\delta$, if
 (\ref{th11}) is satisfied.\\
\ \\
Regarding an implication of the above theorem on the self-similar
blowing up solution to the 3D Euler equations, we have the following
corollary.

\begin{cor}
Let $v\in C([0, T);B^1_{\infty, 1}(\Bbb R^3))$  be a classical
solution to the 3D Euler equations. There exists $\eta
>0$ such that if we have representation for the vorticity $\o=$curl $v$ by
$$\o(x,t)=\frac{1}{T-t} {\bar{\O}}
\left(\frac{x}{(T-t)^{\frac{1}{\a+1}}} \right) \quad \forall
(x,t)\in \Bbb R^3 \times (t_0, T)
$$
for some $t_0 \in (0, T)$, where $\bar{\O}$=curl $\bar{V}$ satisfies
  $\|\bar{\O}\|_{\dot{B}^0_{\infty ,0}} <\eta$, then $\bar{\Omega}=0$, and
  $v\in C([0, T+\delta);B^1_{\infty, 1}(\Bbb R^3))$ for some $\delta
 >0$.
\end{cor}

The following theorem exclude the possibility of a type of
asymptotically self-similar singularity for the 3D Euler equations.
\begin{theorem}
Let $v\in C([0, T);B^1_{\infty, 1}(\Bbb R^3))$  be a classical
solution to the 3D Euler equations.
 Suppose there exist $p_1>0$, $\a > -1$, $\bar{V}\in C^1 (\Bbb R^3)$ such
that $\bar{\Omega}=$curl $\bar{V}\in L^q(\Bbb R^3)$ for all $q\in
(0,p_1)$, and
  \bb\label{th12}
   \lim_{t\nearrow T}
(T-t)\left\|\o(\cdot, t) -\frac{1}{T-t} \bar{\O}
\left(\frac{\cdot}{(T-t)^{\frac{1}{\a+1}}}
\right)\right\|_{\dot{B}^0_{\infty, 1} } =0.
  \ee
 Then, $\bar{\O}=0$, and $v\in C([0, T+\delta);B^1_{\infty, 1}(\Bbb R^3))$ for some $\delta
 >0$.
\end{theorem}
\noindent{\textsf{Remark 1.3}} Although we used the Besov space
$\dot{B}^0_{\infty, 1}(\Bbb R^3)$ for the vorticities in Theorem
1.1 and Theorem 1.2, it would be interesting to see if one could
prove similar results with $\dot{B}^0_{\infty, 1}(\Bbb R^3)$
replaced
by a slightly larger space $L^\infty (\Bbb R^3)$.\\

 \noindent{\textsf{Remark 1.4}} We note that Theorem 1.2 still
 does not exclude the possibility that the vorticity convergence to the
 asymptotically self-similar singularity is weaker than $\dot{B}^0_{\infty, 1}$
 sense. Namely, a self-similar vorticity profile could be approached
 from a local classical solution in the pointwise sense in space,
  or in the $L^p(\Bbb R^3)$ sense for
some $p$
 with $1\leq p\leq\infty$.\\

 \subsection{The Navier-Stokes Equations}
Here we are concerned  on the following 3D Navier-Stokes
equations.
 \[
\mathrm{ (NS)}
 \left\{ \aligned
 &\frac{\partial v}{\partial t} +(v\cdot \nabla )v =-\nabla p  + \Delta v,
 \quad (x,t)\in {\Bbb R^3}\times (0, \infty) \\
 &\textrm{div }\, v =0 , \quad (x,t)\in {\Bbb R^3}\times (0,
 \infty)\\
  &v(x,0)=v_0 (x), \quad x\in \Bbb R^3
  \endaligned
  \right.
  \]

We first state a continuation principle for local in time $L^p
(\Bbb R^3)$ solution of the  Navier-Stokes equations. Below, we
denote the $L^p(\Bbb R^3)$ norm of $f$ by $\|f\|_{L^p}$.
\begin{theorem}
Let $p\in [3, \infty)$, and $v\in C([0, T);L^p (\Bbb R^3 ))$
 be a classical solution to (NS). There exists
a  constant $\eta >0$ depending on $p$ such that if
 \bb\label{th13}
 \inf_{0\leq t<T} (T-t)^{\frac{p-3}{2p}} \|v (t)\|_{L^p} <\eta ,
 \ee
 then, $v\in C([0, T+\delta); L^p(\Bbb R^3 ))$ for some $\delta
 >0$.
\end{theorem}
\noindent{ \textsf{Remark 1.5}} Given $v_0 \in L^p (\Bbb R^3)$
with $p\in [3, \infty)$, the existence and the uniqueness of local
in time classical solution $v\in C([0, T);L^p (\Bbb R^3 ))$  are
established by Kato in \cite{kat};
moreover, the solution is smooth for all $t\in (0, T)$.\\
\ \\
 Similarly to Corollary 1.1 we can reproduce the results of
\cite{nec, tsa} easily under the assumption of additional smallness
condition.
\begin{cor}
Let $p\in [3, \infty)$, and $v\in C([0, T);L^p (\Bbb R^3 ))$
 be a classical solution to (NS). There exists
a  constant $\eta >0$ depending on $p$ such that if
$$
v(x,t)=\frac{1}{\sqrt{T-t}}\bar{V} \left(\frac{x}{\sqrt{T-t}}
\right) \quad \forall (x,t)\in \Bbb R^3 \times (t_0, T)
$$
for some $t_0 \in (0, T)$, where $\bar{V}$ satisfies
  $\|\bar{V}\|_{L^p} <\eta$, then $\bar{V}=0$, and
  $v\in C([0, T+\delta);L^p(\Bbb R^3))$ for some $\delta
 >0$.
\end{cor}

The following theorem for the case $p\in (3, \infty)$ was obtained
by Hou and Li in \cite{hou}. In the next section we present an
alternative proof, which is very simple and elementary compared to
the one given in \cite{hou}.

\begin{theorem} Let $p\in [3, \infty)$, and $v\in C([0, T);L^p (\Bbb R^3 ))$
 be a classical solution to (NS).
 Suppose there exists $\bar{V}\in L^p(\Bbb R^3)$ such
that
  \bb \label{th14}
  \lim_{t\nearrow T}
(T-t)^{\frac{p-3}{2p}}\left\|v(\cdot, t) -
\frac{1}{\sqrt{T-t}}\bar{V} \left(\frac{\cdot}{\sqrt{T-t}}
\right)\right\|_{L^p} =0.
  \ee
 Then, $\bar{V}=0$, and $v\in C([0, T+\delta); L^p (\Bbb R^3 ))$ for some $\delta
 >0$.
\end{theorem}

We now consider a version of localization of Theorem 1.4, in which
we consider the case where the local smooth solution converges to
a self-similar profile in a shrinking ball with the radius
proportional to $\sqrt{T-t}$ as $t\nearrow T$. We denote
$B(z,r)=\{ x\in \Bbb R^3\,|\, |x-z|<r\}$ below.

\begin{theorem}  Let $p\in [3, \infty)$, and $v\in C([0, T);L^p (\Bbb R^3 ))$
 be a classical solution to (NS). Suppose either one of the
 followings hold.
 \begin{itemize}
 \item[(i)] Let $q\in [3, \infty)$. Suppose  there
exists $\bar{V}\in L^p (\Bbb R^3 )$ and $R\in (0, \infty)$ such
that  we have
  \bq \label{th15}
  \lim_{t\nearrow T}\,
(T-t)^{\frac{q-3}{2q}}\sup_{t<\tau <T}\left\|v(\cdot, \tau)
-\frac{1}{\sqrt{T-\tau}} \bar{V} \left(\frac{\cdot
-z}{\sqrt{T-\tau}} \right)\right\|_{L^q(B(z,R\sqrt{T-t}\,))}
=0.\n \\
  \eq
  \item[(ii)]  Let $q\in [2, 3)$. Suppose  there
exists $\bar{V}\in L^p (\Bbb R^3 )$ such that (\ref{th15}) holds
for all $R\in (0, \infty)$.
\end{itemize}
 Then, $\bar{V} =0$, and $v(x,t)$ is H\"{o}lder continuous near
$(z,T)$ in the space and the time variables.
\end{theorem}
\noindent{\textsf{Remark 1.7}}  We note that, in contrast to Theorem
1.4, the range of $q\in [2, 3)$ is also allowed for the possible
convergence of the local classical solution to the self-similar
profile.\\
\ \\
As an immediate corollary of Theorem 1.5(i) we have the following
result, which is a local version of \cite{nec, tsa}.
\begin{cor}
Let $p\in [3, \infty)$, and $v\in C([0, T);L^p (\Bbb R^3 ))$
 be a classical solution to (NS). Suppose there exists $\bar{V}\in L^p (\Bbb R^3 )$ and
  $r>0$ such that
we have the representation,
 $$
 v(x, t)
=\frac{1}{\sqrt{T-t}} \bar{V} \left(\frac{x -z}{\sqrt{T-t}} \right)
\quad \forall (x,t)\in B(z,r)\times (T-r^2, T).
$$
Then, $\bar{V} =0$, and $v(x,t)$ is H\"{o}lder continuous near
$(z,T)$ in the space and the time variables.
\end{cor}
\noindent{\textsf{Remark 1.8}} We note that the above corollary can
be also deduced by a different reasoning from the above, based on
the results of \cite{nec, tsa} combined with simple scaling
argument, which is done in \cite{cha2}.

\section{Proofs for the Euler equations}
 \setcounter{equation}{0}

\noindent{\bf Proof of Theorem 1.1} We start from the following
basic a priori estimates in $B^1_{\infty, 1} (\Bbb R^3)$ and
$\dot{B}^1_{\infty, 1} (\Bbb R^3)$ for the solution $v\in C([0,
T); B^1_{\infty, 1} (\Bbb R^3 ))$ of the Euler equations(see e.g.
\cite{cha3, cha5}), which is a preliminary step to establish the
local existence in $B^1_{\infty, 1}(\Bbb R^3)$. For the
inhomogeneous norm we have
  \bb\label{inhom}
\frac{d}{dt} \|v(t)\|_{{B}^1_{\infty, 1}} \leq C \|\nabla v
(t)\|_{L^\infty} \|v(t)\|_{{B}^1_{\infty, 1}}.
 \ee
 By Gronwall's lemma applied to (\ref{inhom}) we obtain
 \bq\label{bkm}
 \|v(t)\|_{{B}^1_{\infty, 1}} &\leq &\|v_0\|_{{B}^1_{\infty, 1}}
 \exp\left[ C\int_0 ^t \|\nabla v(\tau)\|_{\dot{B}^1_{\infty, 1}}d\tau\right]\n\\
 &\leq&
\|v_0\|_{{B}^1_{\infty, 1}}
 \exp\left[ C\int_0 ^t \|\o (\tau)\|_{\dot{B}^1_{\infty,
 1}}d\tau\right],
 \eq
  where we used the embedding inequality,
 $\dot{B}^0_{\infty, 1}
(\Bbb R^3) \hookrightarrow L^\infty (\Bbb R^3)$.
  On the other hand, for the homogeneous norms we have
 \bb\label{hom}
   \frac{d}{dt} \|v(t)\|_{\dot{B}^1_{\infty, 1}} \leq
C \|\nabla v (t)\|_{L^\infty} \|v(t)\|_{\dot{B}^1_{\infty, 1}}\leq
C \|v(t)\|_{\dot{B}^1_{\infty, 1}} ^2.
 \ee
 The Gronwall lemma applied to (\ref{hom})
 provides us with
 \bb\label{gro} \|v(t)\|_{\dot{B}^1_{\infty, 1}} \leq
\frac{\|v_0\|_{\dot{B}^1_{\infty,
1}}}{1-Ct\|v_0\|_{\dot{B}^1_{\infty, 1}}}.
  \ee
  Translating in time, using the equivalence of
norms, $\|v(t)\|_{\dot{B}^1_{\infty, 1}} \simeq \|\o (t)\bn$, we
have instead of (\ref{gro})
 \bb\label{pro11}
 \|\o(T_1)\bn
\leq \frac{C_1\|\o (t)\|_{\dot{B}^0_{\infty, 1}}}{1-C_2(T_1-t)\|\o
(t) \bn}
 \ee
 for all $T_1 \in (0,T)$ and $t\in [0, T_1)$, where   $C_1, C_2$ are absolute
 constants.
We set $\eta=\frac12 C_2$. For such $\eta$, we  suppose
(\ref{th11}) holds true. Then,  there exists $t_1 \in  [0, T)$
such that $(T-t_1)\|\o (t_1)\bn <\eta$. Fixing $t=t_1$, and
passing $T_1 \nearrow T$ in (\ref{pro11}), we find that
$$\lim\sup_{T_1 \nearrow T} \|\o (T_1 )\bn \leq  2C_1 \|\o( t_1 )\bn
<\infty .$$
  Hence, $\int_0 ^T \|\o (t)\bn dt <\infty$, and  by (\ref{bkm}), we firstly have
  a continuation of local in time solution $v(\cdot, t)$ up to $t=T$, and   $v\in C([0, T];B^1_{\infty, 1} (\Bbb R^3))$.
By the local existence
  theorem applied to the initial data at $T$, we can  continue  further the solution
  $v(\cdot, t)$ until $t=T+\delta $ for some $\delta >0$. $\square$\\
\ \\
\noindent{\bf Proof of Corollary 1.1} We just observe that
$(T-t)\|\o (t)\bn = \|\bar{\O}\bn$ for all $t\in (0, T)$. Hence, our
smallness condition, $\|\bar{\Omega}\bn <\eta $, implies that
$\inf_{0<t<T}(T-t)\|\o (t)\bn <\eta$. Applying Theorem 1.1, we
conclude our proof. $\square$\\
\ \\
\noindent{\bf Proof of Theorem 1.2} We change variables from $(x,t)
\in \Bbb R^3 \times [0,T)$ into $(y,s)\in \Bbb R^3\times [0,
\infty)$ as follows:
$$
y=\frac{x}{(T-t)^{\frac{1}{\a+1}}}, \quad s=\frac{1}{\a+1} \log
\left( \frac{T}{T-t}\right).
$$
We note that this type of  introduction of similarity variables
was previously used in \cite{gig2} in the context of nonlinear
heat equation. Based on this change of variables, we transform
$(v,p)\mapsto (V, P)$ according to
 \bb
 v(x,t)=\frac{1}{(T-t)^\frac{\a}{\a+1}} V(y,s ), \quad
 p(x,t)=\frac{1}{(T-t)^\frac{2\a}{\a+1}} P(y,s ).
 \ee
 Substituting $(v,p)$ into the $(E)$ we obtain the equivalent evolution equation for
 $(V,P)$,
 $$
 (E_1) \left\{ \aligned
  &  V_s +\a V +(y \cdot \nabla)V +(\a+1) (V\cdot \nabla )V =-\nabla
  P,\\
 & \mathrm{div}\, V=0,\\
 & V(y,0)=V_0 (y)=T^{\frac{\a}{\a +1}} v_0 (T^{\frac{\a}{\a +1}}y).
  \endaligned \right.
 $$
 In terms of $V$ the condition (\ref{th12}) is translated
into
 \bb\label{econv1}
 \lim_{s\to \infty}\|\O(\cdot ,s) -\bar{\O}(\cdot )\bn=0,
 \ee
 where we set $\O =$ curl $V$.
 Combining this with  the embedding, $\b (\Bbb R^3) \hookrightarrow L^\infty (\Bbb R^3)$
 and the fact that the Calderon-Zygmund singular integral operator maps $\b (\Bbb R^3)$ into
 itself boundedly, we obtain
 \bb\label{econv2}
  \lim_{s\to \infty}\|\nabla V(\cdot ,s) -\nabla \bar{V}(\cdot )\|_{L^\infty}=0.
  \ee
Similarly to \cite{hou}, we consider the scalar test function $\xi
\in C^1_0 (0,1)$ with $\int_0 ^1\xi (s)ds\neq 0$,  and the vector
test function $\phi =(\phi_1 , \phi_2, \phi_3 )\in C_0^1 (\Bbb
 R^3)$ with div $\phi=0$. We multiply the first equation of $(E_1)$
 in the dot product by $\xi
 (s-n)\phi (y)$, and integrate it over $\Bbb R^3\times [n, n+1]$,
 and then we  integrate by part for the terms including the
 time derivative and the for the pressure term to obtain
 \bqn
 &&-\int_0^{1}\int_{\Bbb R^3} \xi _s(s) \phi(y)\cdot V(y,s+n)
 dyds\\
 &&+\int_0 ^{1}\int_{\Bbb R^3}\xi (s)\phi(y) \cdot[\a V +(y \cdot \nabla)V +(\a+1)
 (V\cdot \nabla )V](y,s+n)  dyds=0.
 \eqn
  Passing to the limit $n\to \infty$ in this equation, using the fact
  (\ref{econv2}), $\int_0 ^1\xi _s(s)ds=0$ and $\int_0 ^1\xi
  (s)ds\neq 0$,
  we  find that $\bar{V}\in C^1 (\Bbb R^3)$ satisfies
$$
\int_{\Bbb R^3} [\a \bar{V} +(y \cdot \nabla)\bar{V} +(\a+1)
(\bar{V}\cdot \nabla
 )\bar{V}]\cdot \phi dy=0
 $$
 for all $\phi \in C_0^1 (\Bbb
 R^3)$ with div $\phi=0$. Hence, there exists a scalar function
 $\bar{P}$ such that
 \bb\label{eleray1}
\a  \bar{V}+ (y\cdot\nabla )\bar{ V}
 +(\a+1)(\bar{V}\cdot \nabla )\bar{V} =-\nabla  \bar{P}.
 \ee
 On the other hand, we can pass $s\to \infty$ directly in the
 second equation of $(E_1)$ to have
 \bb\label{eleray2}
 \mathrm{div}\, \bar{V}=0.
 \ee
  Since $\bar{V}$ is a classical solution of (\ref{eleray1})-(\ref{eleray2}),
  and curl $\bar{V}=\bar{\Omega}$
  satisfy the conditions of Theorem 1.1 of \cite{cha1} by our hypothesis, we can deduce
  that $\bar{\Omega}=0$ by that theorem.
 Hence,  (\ref{econv1}) implies that $\lim_{s\to \infty} \|\O
 (s)\bn=0$.
 Thus, for $\eta >0$ given in Theorem 1.1,
 there exists $s_1>0$ such that $\|
 \O (s_1 )\bn <\eta $.
Let us set $t_1=T[1-e^{(\a +1)s_1} ]$. Going back to the original
physical variables, we have
 $$
  (T-t_1)\|\o (t_1)\bn <\eta .
 $$
 Applying Theorem 1.1, we conclude the proof. $\square$

\section{Proofs for the Navier-Stokes equations}
 \setcounter{equation}{0}

\noindent{\bf Proof of Theorem 1.3} In the case $p=3$, the
conclusion of Theorem 1.3 follows immediately from the result of
the small data global existence  in $L^3 (\Bbb R^3)$, proved by
Kato(\cite{kat}). Below we concentrate on the case $p\in (3,
\infty)$. We recall the following result for the blow-up rate
estimate essentially obtained by Leray(pp.227,\cite{ler}): Suppose
$v\in C([0,T); L^p (\Bbb R^3))$, $p\in (3,\infty)$,  is a local in
time classical solution to (NS). Then, we have
 \bb\label{leray}
 \quad \lim\sup_{t\nearrow T}
\|v(t)\|_{L^p} =\infty \Rightarrow \|v(t)\|_{L^p} \geq \frac{K}{(T
-t)^{\frac{p-3}{2p}}} \ee
 for all $t\in [0, T)$ with a constant $K=K(p)$ independent of $T$ and $t$.

Let $\kappa$ be the supremum of the constant $K$ in (\ref{leray}).
Then, choosing $\eta =\kappa/2$, and
 taking  the contraposition of the statement, we
deduce the following: If there exists $t_1 \in [0,T)$ such that
$(T-t_1)^{\frac{p-3}{2p}} \|v(t_1)\|_{L^p}<\eta$, then
$\lim\sup_{t\nearrow T}\|v(t)\|_{L^p} <\infty$. Now the condition
(\ref{th13}) implies that there exists really such $t_1$, and
hence we have $\sup_{0<t< T} \|v(t)\|_{L^p} <\infty. $ Applying
the local existence result in $L^p (\Bbb R^3)$ due to
Kato(\cite{kat}) with the initial data $v(\cdot,T-\vare) \in L^p
(\Bbb R^3)$ for sufficiently small $\vare
>0$, we can  extend the solution to be $v \in C([0, T+\delta )
;L^p (\Bbb R^3))$ for some $\delta >0$. $\square$\\
\ \\
\noindent{\bf Proof of Corollary 1.2} Similarly to the proof of
Corollary 1.1, we just observe that $(T-t)^{\frac{p-3}{2p}}\|v
(t)\|_{L^p} = \|\bar{V}\|_{L^p}$ for all $t\in (0, T)$. Hence, our
smallness condition, $\|\bar{V}\|_{L^p} <\eta $, implies that
$\inf_{0<t<T}(T-t)^{\frac{p-3}{2p}}\|v(t)\|_{L^p} <\eta$. Applying
Theorem 1.3, we
conclude our proof. $\square$\\
\ \\
\noindent{\bf Proof of Theorem 1.4}
 Similarly to the proof of Theorem 1.2 we change variables from
$(x,t) \in \Bbb R^3 \times [0,T)$ into $(y,s)\in \Bbb R^3\times [0,
\infty)$ as follows.
 \bb\label{simv1}
   y=\frac{x}{\sqrt{T-t}}, \quad s=\frac{1}{2}
\log \left( \frac{T}{T-t}\right).
 \ee
  Based on this change of
variables, we transform $(v,p)\mapsto (V, P)$ according to
 \bb\label{simv2}
 v(x,t)=\frac{1}{\sqrt{T-t}} V(y,s ), \quad
 p(x,t)=\frac{1}{T-t} P(y,s ).
 \ee
 Substituting $(v,p)$ into the $(E)$ we obtain the following equivalent evolution equations for
 $(V,P)$,
 $$
 (NS_1) \left\{ \aligned
  & V_s + V +(y \cdot \nabla)V +2 (V\cdot \nabla )V =-\nabla
  P +2\Delta V,\\
 & \mathrm{div}\, V=0,\\
 & V(y,0)=V_0 (y)=\sqrt{T} v_0 (\sqrt{T} y).
  \endaligned \right.
 $$
 In terms of $V$ the condition (\ref{th14}) is translated
into
 \bb\label{nsconv}
 \lim_{s\to \infty}\|V(\cdot ,s) -\bar{V}(\cdot )\|_{L^p}=0.
 \ee
 Using this convergence, we can pass to the limit $s\to \infty$ in the weak formulation of
$(NS_1)$, as is done in \cite{hou},  which is also an obvious
modification of the step described in the proof of Theorem 1.2.
Thus, we find that
  $\bar{V}$ is a weak solution of the following stationary Leray system(\cite{ler}),
 \bb\label{ler}
 \left\{\aligned & \bar{V}+ (y\cdot\nabla )\bar{ V}
 +2(\bar{V}\cdot \nabla )\bar{V} =-\nabla  \bar{P} +2\Delta \bar{V},
 \\
  &\textrm{div }\bar{ V} =0.
  \endaligned\right.
 \ee
Since $\bar{V}\in L^p (\Bbb R^3)$ by hypothesis,  thanks to the
results in \cite{nec,tsa}(specifically, we use the result of
\cite{nec} for
  $p=3$, while we use result of \cite{tsa} for
  $p>3$),  we can deduce that $\bar{V}=0$.
 Hence,  (\ref{nsconv}) implies that $\lim_{s\to \infty}
 \|V(s)\|_{L^p}=0$ and, for  $\eta >0$ given in Theorem 1.3, there exists $s_1>0$ such
that
 \bb\label{small}
 \|V(s_1 )\|_{L^p} <\eta .
  \ee
   Let us set $t_1=T[1-e^{2 s_1} ]$. Going
back to the original physical variables,  we can rewrite
(\ref{small}) as
 $$
  (T-t_1)^{\frac{p-3}{2p}}\|v(t_1)\|_{L^p} <\eta .
 $$
 Applying Theorem 1.3, we conclude the proof. $\square$\\
 \ \\
 In order to prove Theorem 1.5 we recall the following
 recent result due to Gustafson, Kang and Tsai(Theorem 1.1,\cite{gus}),
 which we state a part of the theorem
in an easily applicable form for our purpose.  For the statement
of the theorem we recall that a suitable weak solution of the
Navier-Stokes equations is a pair $(v,p)$, satisfying the
equations in (NS) in the sense of distribution, and satisfying the
generalized energy inequality. For more precise definition and its
global in time construction we refer \cite{caf}(see also
\cite{lin} for a refined definition, regarding the integrability
of pressure).

 \begin{theorem} Let $q\in (3/2, \infty)$. Suppose $v$ is a suitable
 weak solution of (NS) in a cylinder, say  $Q=B(z,r_1)\times (t-r_1^2, T)$ for some $r_1>0$.
 Then, there exists a constant
 $\eta =\eta(q)>$ such that if
 \bb\label{kang}
 \lim\sup_{r\searrow 0}\left\{r^{\frac{q-3}{q}} \mathrm{ess}\sup_{t-r^2 <\tau <
 t}\|v(\cdot,\tau)\|_{L^q (B(z, r))}\right\} \leq \eta,
 \ee
 then   $v$ is H\"{o}lder continuous both in space and time
 variables near $(z,t)$.
 \end{theorem}
 \ \\
\noindent{\bf Proof of Theorem 1.5} We first claim that in the
case $q\in [3, \infty)$ the condition (\ref{th15}) for some $R\in
(0, \infty)$ is equivalent to
 \bq
\label{convv14}
  \lim_{t\nearrow T}\,
(T-t)^{\frac{q-3}{2q}}\sup_{t<\tau <T}\left\|v(\cdot, \tau)
-\frac{1}{\sqrt{T-\tau}} \bar{V} \left(\frac{\cdot
-z}{\sqrt{T-\tau}} \right)\right\|_{L^q(B(z,R\sqrt{T-t}\,))}
=0\n \\
  \eq
{\em for all } $R\in (0, \infty)$. Indeed, suppose $R_1 \leq R_2$ be
given, then setting
$$ f(x,\tau):=v(x, \tau)
-\frac{1}{\sqrt{T-\tau}} \bar{V} \left(\frac{x-z}{\sqrt{T-\tau}}
\right), $$
 we have an obvious inequality
 \bq
\label{convv14a}
  &&\lim\sup_{t\nearrow T}\,
(T-t)^{\frac{q-3}{2q}}\sup_{t<\tau <T}\left\|f(\cdot
,\tau)\right\|_{L^q(B(z,R_1\sqrt{T-t}\,))}\n \\
&&\qquad\qquad\leq \lim\sup_{t\nearrow
T}\,(T-t)^{\frac{q-3}{2q}}\sup_{t<\tau <T}\left\|f(\cdot
,\tau)\right\|_{L^q(B(z,R_2\sqrt{T-t}\,))}.
 \eq
On the other hand,  let us suppose $R_2 \geq R_1$. Then,  we observe
the inclusion relation,
 \bb \label{convv14b}
  B(z,R_2
\sqrt{T-t_2}\,) \subset B(z,R_1 \sqrt{T-t_1}\,)\quad
\mbox{if}\quad t_2 \geq T-\left(\frac{R_1}{R_2}\right)^2 (T-t_1).
 \ee
 Since
 $T-\left(\frac{R_1}{R_2}\right)^2 (T-t_1)\geq
 t_1$ for $R_2 \geq R_1$,
 we also have  $t_2\geq t_1$.
  Hence, given $t_1\in (0, T)$ and  $t_2\in (t_1 ,T)$ satisfying (\ref{convv14b}), we have
 \bq\label{convv14c}
&&(T-t_2)^{\frac{q-3}{2q}}\sup_{t_2<\tau <T}\left\|f(\cdot ,\tau)\right\|_{L^q(B(z,R_2\sqrt{T-t_2}\,))}\n\\
 &&\qquad\qquad \leq (T-t_1)^{\frac{q-3}{2q}}\sup_{t_1<\tau <T}
 \left\|f(\cdot ,\tau)\right\|_{L^q(B(z,R_1\sqrt{T-t_1}\,))}
 \eq
if $q\geq 3$. Taking $\lim\sup_{t_2 \nearrow T} $ in the left hand
side of (\ref{convv14c}), and then taking $\lim\sup_{t_1 \nearrow
T}$ in the right hand side of it, we obtain
 \bq
\label{convv14d}
 && \lim\sup_{t\nearrow T}\,
(T-t)^{\frac{q-3}{2q}}\sup_{t<\tau <T}\left\|f(\cdot ,\tau)\right\|_{L^q(B(z,R_2\sqrt{T-t}\,))}\n\\
&&\qquad\qquad\leq \lim\sup_{t\nearrow
T}\,(T-t)^{\frac{q-3}{2q}}\sup_{t<\tau <T}\left\|f(\cdot
,\tau)\right\|_{L^q(B(z,R_1\sqrt{T-t}\,))}.
 \eq
 Combining (\ref{convv14a}) and (\ref{convv14d}), we have
\bqn
 && \lim\sup_{t\nearrow T}\,
(T-t)^{\frac{q-3}{2q}}\sup_{t<\tau <T}\left\|f(\cdot ,\tau)\right\|_{L^q(B(z,R_2\sqrt{T-t}\,))}\\
&&\qquad=\lim\sup_{t\nearrow T}\,(T-t)^{\frac{q-3}{2q}}\sup_{t<\tau
<T}\left\|f(\cdot ,\tau)\right\|_{L^q(B(z,R_1\sqrt{T-t}\,))}
 \eqn
for all $0<R_1<R_2<\infty$ if $q\geq 3$, and our claim is proved.
Hence, for all $q\in [2, \infty)$ either the condition (i), or
condition (ii) of Theorem 1.5 implies that (\ref{th15}) holds
  for all $R>0$.
As previously we change variables from $(x,t) \in B(z, R\sqrt{T-t}\,
) \times [0, T) $ into $(y,s)\in
 B(0, R) \times [0, \infty)$ as follows.
 \bb
   y=\frac{x-z}{\sqrt{T-t}}, \quad s=\frac{1}{2}
\log \left( \frac{T}{T-t}\right).
 \ee Based on this change of
variables, we transform $v\mapsto V$ according to
 \bb
 v(x,t)=\frac{1}{\sqrt{T-t}} V(y,s ).
 \ee
Then, (\ref{th15}) is written as
 \bb\label{convv15}
 \lim_{s\to \infty} \|V(\cdot, s)-\bar{V}\|_{L^q(B(0, R))}=0,
 \ee
 which holds for all $R\in (0, \infty)$ by the  above claim.
 Following exactly the same procedure as in the proof of Theorem
 1.4, we can conclude that $\bar{V}$ is a weak solution
 of the Leray system (\ref{ler}). We note here that $L^q_{loc} (\Bbb R^3)$
 convergence with $q\in [2,
 \infty)$ is enough to show that $\bar{V}$ a weak solution of the
 Leray system.
On the other hand, by hypothesis $\bar{V } \in L^p (\Bbb R^3 )$
 with $p\in [3, \infty)$, hence $\bar{V}=0$ by the results of \cite{nec,tsa}.
Therefore, (\ref{th15}) is reduced to
 \bb \label{pro15}
  \lim_{t\nearrow T}\,
(T-t)^{\frac{q-3}{2q}}\sup_{t<\tau<T}\left\|v(\cdot,
\tau)\right\|_{L^q(B(z,R\sqrt{T-t}\,))} =0
  \ee
  for all $R\in (0, \infty)$.
 We set $R=1$ and $\sqrt{T-t}=r$ in
  (\ref{pro15}), which becomes
 \bb \label{pro16}
\lim_{r\searrow 0}\left\{
r^{\frac{q-3}{q}}\sup_{T-r^2<\tau<T}\left\|v(\cdot,
\tau)\right\|_{L^q(B(z,r))}\right\} =0.
  \ee
 Hence, the conclusion follows from Theorem 3.1. $\square$\\
\[ \mbox{\bf Acknowledgements}\]
The author would like to thank to Professor Tai-Peng  Tsai for
helpful discussion.


\begin{thebibliography}{1}
\bibitem{bea}J. T. Beale, T. Kato and A. Majda,  {\it Remarks on the
breakdown of smooth solutions for the 3-D Euler equations}, Comm.
Math. Phys., {\bf 94},  (1984), pp. 61-66.
\bibitem{caf} L. Caffarelli, R. Kohn and L. Nirenberg,
{\it Partial regularity of suitable weak solutions of the
Navier-Stokes equations,} Comm. Pure Appl. Math., {\bf 35},
(1982), pp.771-831.
\bibitem{cha1} D. Chae, {\it Nonexistence of
self-similar singularities for the 3D incompressible Euler
equations,} arXiv-preprint, math.AP/0601060.
\bibitem{cha2} D. Chae, {\it A note on  `Nonexistence of  self-similar singularities
for the 3D incompressible Euler equations'}, arXiv-preprint,
math.AP/0601661.
 %\bibitem{chaa}D. Chae, {\it Remarks on the blow-up criterion of the 3D Euler
 %equations,} Nonlinearity, {\bf 18}, (2005), pp. 1021-1029.
 \bibitem{cha3}D. Chae, {\it  Local Existence and Blow-up Criterion for the
Euler Equations in the Besov Spaces}, Asymp. Anal., {\bf 38}, no.
3-4, (2004), pp. 339-358.
\bibitem{cha4} D. Chae, {\it On the Euler Equations in the Critical
 Triebel-Lizorkin Spaces},
  Arch. Rat. Mech. Anal., {\bf 170}, no. 3, (2003), pp.185-210.
 \bibitem{cha5} D. Chae, {\it Remarks on the blow-up of the Euler equations
 and the related equations},
   Comm. Math. Phys., {\bf 245}, no. 3, (2003), 539-550.
 \bibitem{che}J. Y. Chemin, {\it Perfect incompressible fluids,}
 Clarendon Press, Oxford, (1998).
\bibitem{con1}P. Constantin, {\it Geometric Statistics in
Turbulence}, SIAM Rev.,{\bf 36}, (1994), pp. 73-98.
\bibitem{con2}P. Constantin,  {\it A few results and open problems
regarding incompressible fluids,} Notices Amer. Math. Soc. {\bf
42}, no. 6, (1995), pp. 658-663.
\bibitem{con3}P. Constantin, C. Fefferman and A. Majda,
 {\it Geometric constraints on potential singularity formulation in the
 3-D Euler equations}, Comm. P.D.E, {\bf 21}, (3-4), (1996), pp.
 559-571.
 \bibitem{con4}P. Constantin and C. Foias, {\it Navier-Stokes Equations,}
 Chicago Lectures in Mathematics Series, Univ. Chicago Press (1988).
 \bibitem{esc} L. Escauriaza, G. Seregin and V.
Sverak, {\it $L^{3,\infty}$-solutions of Navier-Stokes equations
and backward uniqueness,} Rus. Math. Surveys, {\bf 58}, (2003),
pp. 211-250.
%\bibitem{gal}G. Galdi, {\it An Introduction to the mathematical
%theory of Navier-Stokes equations, I,II}, Springer-Verlag, (1994).
 %\bibitem{gig1}Y. Giga, {\it Solutions for semilinear parabolic
 %equations in $L^p$ and regularity of weak solutions of the
 %Navier-Stokes system,} J. Diff. Eqns., {\bf 62}, (1986), pp.
 %186-212.
 \bibitem{gig2}Y. Giga and R. V. Kohn, {\it Asymptotically
 Self-Similar Blow-up of Semilinear Heat Equations,} Comm. Pure
 Appl. Math., {\bf 38}, (1985), pp. 297-319.
 \bibitem{gre}J. M. Greene and R. B. Pelz, {\it Stability of postulated,
self-similar, hydrodynamic blowup solutions}, Phys. Rev. E, {\bf
62}, no. 6, pp. 7982-7986.
\bibitem{hou}T. Y. Hou and R. Li, {\it Nonexistence of Local Self-Similar Blow-up for the 3D
Incompressible Navier-Stokes Equations,} arXiv-preprint,
math.AP/0603126.
\bibitem{kat} T. Kato, {\it Strong $L^p$ solutions of the
Navier-Stokes equations in ${\mathbb R}^m$ with applications to
weak solutions,} Math. Zeit., {\bf 187}, (1984), pp. 471-480.
\bibitem{ker}R. Kerr, {\it Evidence for a singularity of the
3-dimensional, incompressible Euler equations,} Phys. Fluids A, {\bf
5}, (1993), pp. 1725-1746.
\bibitem{lad}O. A. Ladyzenskaya, {\it The mathematical theory of
viscous incompressible flow,} Gordon and Breach, (1969).
\bibitem{ler}J. Leray, {\it Essai sur le mouvement d'un fluide
visqueux emplissant l'espace,} Acta Math. {\bf 63} (1934), pp.
193-248.
\bibitem{lin}F. Lin, {\it A
new proof of the Caffarelli-Kohn-Nirenberg theorem,} Comm. Pure
Appl. Math. {\bf 51},  no. 3, (1998), pp. 241-257.
\bibitem{maj1}A. Majda, {\it Vorticity and the mathematical theory
of incompressible fluid flow,} Comm. Pure Appl. Math., {\bf 39},
(1986), pp. 187-220.
\bibitem{maj2}A. Majda and A. Bertozzi, {\it Vorticity and
Incompressible Flow,} Cambridge Univ. Press. (2002).
\bibitem{mil}J. R. Miller, M. O'Leary and M. Schonbek, {\it
Nonexistence of singular pseudo-self-similar solutions of the
Navier-Stokes system,} Math. Ann. {\bf 319}, (2001), no. 4, pp.
809-815.
\bibitem{nec}J. Necas, M.  Ruzicka and V. Sverak, {\it On Leray's
self-similar solutions of the
Navier-Stokes equations,} Acta Math. {\bf 176}, no. 2, (1996), pp.
283-294.
\bibitem{oky}T. Okyama, {\it Interior regularity of weak solutions to the
Navier-Stokes equation,} Proc. Japan Acad. {\bf 36}, (1960), pp.
%273-277.
 \bibitem{pak}H. C. Pak, {\it  Existence of solution for the Euler equations in
a critical Besov space $B^1_{\infty, 1} (\Bbb R^n)$,}
    Comm. P.D.E., {\bf 29}, pp. 1149-1166.
\bibitem{pel}R. Pelz, {\it Symmetry and hydrodynamic blow-up
problem}, J. Fluid Mech., {\bf 444}, pp. 299-320.
\bibitem{pro}G. Prodi, {\em Un Teorema di Unici\'{a} per le Equazionni
di Navier-Stokes,} Ann. Mat. Pura Appl., {\bf 48}, no. 4, (1959),
pp. 1773-182.
\bibitem{ser}J. Serrin, {\it On the Interior Regularity of Weak Solutions of
the Navier-Stokes Equations}, Arch. Rat. Mech. Anal., {\bf 9},
(1962), pp.187-191.
\bibitem{gus}S. Gustafson, K. Kang and T-P. Tsai, {\it Regularity
criteria for suitable weak solutions of the Navier-Stokes equations
near the boundary,} to appear in J. Diff. Eqns, (available online).
\bibitem{tem}R. Temam, {\it Navier-Stokes equations,} 2nd ed.,
 North-Holland, Amsterdam,
(1986).
\bibitem{tri}H. Triebel,  {\it Theory of Function Spaces},
Birk\"{a}user Verlag, Boston, (1983).
\bibitem{tsa}T-P. Tsai, {\it  On Leray's self-similar solutions of the Navier-Stokes
equations satisfying local energy estimates,} Arch. Rat. Mech.
Anal., {\bf 143}, no. 1, (1998), pp. 29-51.
\bibitem{vis}M. Vishik, {\it Hydrodynamics in Besov spaces,}
Arch. Rat. Mech. Anal, {\bf 145}, (1998), pp. 197-214.
  \end{thebibliography}
\end{document}